\documentstyle{article}
\title{Graph and Union of Graphs Compositions}
\author{W. Bajguz\\University of Bialystok}

\newtheorem{definition}{Definition}
\newtheorem{theorem}{Theorem}
\newtheorem{conclusion}{Conclusion}
\newtheorem{lemma}{Lemma}

\begin{document}
\maketitle

\begin{abstract}
The graph compositions` notion was introduced by A. Knopfmacher and
M. E. Mays \cite{KM}. In this note we add to these a new
construction of tree-like graphs where nodes are graphs themselves.
The first examples of these tree-like compositions, a corresponding
theorem and resulting conclusions are provided.
\end{abstract}

\small{\noindent AMS Classification Numbers: 05A18, 05C05, 05C30,
11B37.

\vspace{0.2cm}

\noindent Key Words: graph compositions, union of graphs

\vspace{0.2cm}

\noindent presentation  (January $2006$) at the Gian-Carlo Rota
Polish Seminar\\
$http://ii.uwb.edu.pl/akk/index.html$}

\section{Introduction}

Graphs considered here are finite, undirected and labeled graphs,
with no loops or multiple edges. The edge between $v_1,v_2$ is
$(v_1,v_2)$. The set of vertices of graph $G$ is denoted by $V(G)$
and the set of edges - by $E(G)$.

Let $G$ be a graph and $S\subseteq G$ be a subgraph of $G$. We say
that $S$ is a \textbf{maximal subgraph} of $G$ iff
$E(S)=\{(x,y)\in E(G): x,y\in V(S)\}$. Once we have this notion we
can introduce definition of the composition (equivalent to that
from \cite{KM}):

\begin{definition}
The \textbf{composition} of graph $G$ is a partition
$\{V_1,V_2,...,V_n\}$ of the set $V(G)$, such that each maximal
subgraph of $G$ induced by $V_i$ is connected for $1\leq i\leq n$.
\end{definition}

Let \boldmath{$C(G)$} denotes the number of compositions of the
graph $G$. In particularly simple cases, the number of composition
may be counted immediately. \\
For example, let $E_0$ be a graph with empty set of edges. Then $C(E_0)=1$.\\
For   cycle graph $C_4$ with 4 vertices and 4 edges (tetragon) $C(C_4)=12$.\\
Next - the graphs such as  point, interval and triangle, are
special cases of complete graphs $K_n$ with $n$ vertices
($n=1,2,3$). Since $(x,y)\in E(K_n)$ for any two vertices $x,y\in
V(K_n)$, then the number of compositions of $K_n$ is equal to the
number of partitions of a set with $n$ elements, i.e to the $n$-th
Bell number $B(n)$ (\cite{KM} theorem 2):
\begin{equation}
C(K_n)=B(n) \label{K_n}
\end{equation}

\section{Compositions of the union of graphs}

Let us consider at first the union $G_1\cup G_2$ of two graphs
which are disconnected or are connected by one common vertex.
Since each composition of graph $G_1\cup G_2$ is in one to one
correspondence with a pair of compositions of $G_1$ and $G_2$, we
have (\cite{KM} theorem 3):

\begin{theorem} Let $G_1,G_2$ be graphs with no common edges and at
most one common vertex. Then
\begin{equation}
 C(G_1\cup G_2)=C(G_1)\cup C(G_2)
 \label{vertex}
\end{equation}
\end{theorem}

The following formula was independently proved in \cite{KM} (theorem
4), but it immediately follows from above theorem and the equation
$C((x_1,x_2))=2$.

\begin{conclusion}
Let $G_1,G_2$ be a disconnected graphs and $x_1\in V(G_1)$, $x_2\in
V(G_2)$. Then
\begin{equation}
C(G_1\cup (x_1,x_2)\cup G_2)=2\cdot C(G_1)\cdot C(G_2)
\label{long_vertex}
\end{equation}
\end{conclusion}

Consider now another example of applying the above. In \cite{KM}
was shown, that for tree $T_n$ with $n$ edges ($n-1$ vertices),
the number of compositions is $2^n$. As any  tree may be
constructed inductively from one edge by adding succeeding edges,
the proof of formula is by induction. By analogous construction we
can build a tree $T_n(G)$ of $n$ copies of any graph $G$, where in
succeeding steps we connect copy of $G$ with constructed tree by
common vertex. Then $C(T_n(G))=(C(G))^n$.

For $A\subseteq V(G)$ let $C^+(G,A)$ denotes the number of
compositions of a graph $G$ such that set $A$ is in one part of
composition while $C^-(G,A)$ denotes the number of compositions of
the graph $G$ such that each element of a set $A$ is in another
part of a composition. Then we state.

\begin{lemma}
Let $G,H$ be a graphs, $V(G)\cap V(G)=\{x,y\}$ and $(x,y)\in
E(G)\cap E(H)$. Then $$C(G\cup H)=C^+(G,\{x,y\})\cdot
C^+(H,\{x,y\})+C^-(G,\{x,y\})\cdot C^-(H,\{x,y\})$$
\end{lemma}
\textsc{Proof}. In order to prove this equation, it is sufficient
to observe that\\
$C^+(G\cup H,\{x,y\})=C^+(G,\{x,y\})\cdot C^+(H,\{x,y\})$
and\\
$C^-(G\cup H,\{x,y\})=C^-(G,\{x,y\})\cdot C^-(H,\{x,y\})$.\\
Hence the above equation follows by \\
 $C(G\cup H)=C^+(G\cup H,\{x,y\})+C^-(G\cup H,\{x,y\})$. $\Box$

\begin{definition}
Let $G$ be a graph and $k=(x,y)\in E(G)$. Then $G/k$ denotes graph
obtained from $G$ by removing k from edges, identifying vertices
$x,y$ and identifying edges $(x,z)$ with $(y,z)$ for vertices $z\in
V(G)$ such that $(x,z),(y,z)\in E(G)$.
\end{definition}

With this notion applied all together with Lemma 1. we arrive at
the following theorem.

\begin{theorem}
Let $G,H$ be graphs with exactly two common vertices $x,y$ and one
common edge $(x,y)$. Then
\begin{eqnarray}
C(G\cup H)=C(G)\cdot C(H)+2\cdot C(G/k)\cdot C(H/k)- \label{edge}\\
-C(G)\cdot C(H/k)-C(G/k)\cdot C(H) \nonumber
\end{eqnarray}
\end{theorem}
\textsc{Proof}. The thesis  follows from equations for $G$:\\
$C(G)=C^+(G,\{x,y\}))+C^-(G,\{x,y\}))$, $C^+(G,\{x,y\}))=C(G/k)$\\
and analogical equations for $H$. $\Box$

Theorems above make it possible (in some cases) to count the
number of compositions of graphs constructed inductively just by
attaching  succeeding graphs with common vertex or common edge -
similarly like in the construction of trees.

\section{Trees of graphs}

In this section  we use theorems the preceding section in some
special cases, when $n$ copies of graph $G$ are connected into a
tree-like structure i.e.  nodes of a tree are graphs themselves.

\begin{definition}
The V-tree (respectively E-tree) $T$ of graphs $G_1,G_2,...,G_n$ is
a graph $T=T_n$ constructed as following:
\begin{enumerate}
\item $T_1=G_1$,
\item if V-tree (E-tree) $T_k$ is defined for some $k<n$, then $T_{k+1}$ is obtained
from graphs $T_k$ and $G_{k+1}$ by identifying some vertex (edge) in
$T_k$ and $G_{k+1}$.
\end{enumerate}
\end{definition}

Immediately from equation \ref{vertex} the simple conclusion
follows.
\begin{conclusion}
If $T$ is a V-tree of graphs $G_1,G_2,...,G_n$, then
\begin{equation}
C(T)=\prod_{i=1}^{n}(C(G_i))
\end{equation}
\end{conclusion}

E-trees are more complicated, nevertheless in this case we can
obtain an interesting insight too into the matters via the Theorem
3. and resulting consequences.
\begin{theorem}[technical]
Let $G$ be a graph and $T$ be an E-tree of $n$ copies of graph $G$.
Let $T_1,T_2,...,T_n$ be the sequence of E-trees used to
construction $T$ and $k_1,k_2,...,k_{n-1}$ be a sequence of edges in
$T$, such that $T_{r+1}$ is union of graphs $T_r$ and copy of $G$
connected by edge $k_r$ for $1\leq r<n$. Then
$$
\begin{array}{l}
C(T_{r+1})=(C(G)-C(G/k_r))\cdot C(T_r)+\\
\mbox{\hskip 30mm}+(2\cdot C(G/k_r)-C(G))\cdot
C(T_r/k_r)\\
C(T_{r+1}/k_{r+1})=\\
=\left\{
 \begin{array}{lc}
(2\cdot
 C(G/k_{r+1}/k_r)-C(G/k_{r+1}))\cdot
 C(T_r/k_r)\\
 +(C(G/k_{r+1})-C(G/k_{r+1}/k_r))\cdot C(T_r)
  & k_{r+1}\neq k_r\\
  \\
 C(G/k_r)\cdot C(T_r/k_r) & k_{r+1}=k_r
 \end{array}
 \right.
\end{array}
$$
for $1\leq r<n$.
\end{theorem}
\textsc{Proof}. The thesis follows by \ref{edge}. $\Box$

The above equation may be simplified for "regular" trees.

\begin{conclusion}
Let $K_n$ be the complete graph on $n>2$ vertices and let be given a
sequence $T_1,T_2,...,T_m$ of E-trees used to construction of E-tree
$T=T_m$ of $m$ copies of $K_n$, such that any three different copies
of $K_n$ have no common edge in $T$. Then for $r<m$ holds
\begin{eqnarray}
C(T_{r+1})=(B_n-2\cdot B_{n-1}+2\cdot B_{n-2})\cdot C(T_r)+\label{tree}\\
+((B_{n-1})^2-B_n\cdot B_{n-2})\cdot C(T_{r-1})\nonumber
\end{eqnarray}
\end{conclusion}
\textsc{Proof}. Let $k_1,k_2,...,k_m-1$ be the sequence of common
edges in construction of sequence $T_1,T_2,...,T_m$. At first
observe that for any edge $k\in E(K_n)$, $K_n/k$ is complete graph
on $n-1$ vertices.
Therefore by theorem 3 and equation \ref{K_n}\\
$C(T_{r+1})=(B_n-B_{n-1})C(T_r)+(2B_{n-1}-B_n)C(T_r/k_r)$,\\
$C(T_{r+1}/k_{r+1})=(2B_{n-2}-B_{n-1})C(T_r/k_r)+(B_{n-1}-B_{n-2})C(T_r)$.\\
Moreover $T_r/k_r=T_{r-1}\cup K_n/k_r=T_{r-1}\cup K_{n-1}$ and
graphs $T_{r-1},K_{n-1}$ have common exactly one edge $k_{r-1}$ for
$r>1$. From first equation we can compute $C(T_(r-1)\cup K_{n-1})$
and use it in second equation to obtain the thesis. $\Box$

Similarly , although more easily, the following conclusion becomes
apparent.
\begin{conclusion}
Let $K_n$ be the complete graph on $n>2$ vertices and let be given a
sequence $T_1,T_2,...,T_m$ of E-trees used to construction E-tree
$T$ of $m$ copies of $K_n$ with property: there is an edge $k\in
E(T)$ such that every two different copies of $K_n$ have common edge
$k$. Then for $r<m$ holds
\begin{equation}
C(T_{r+1})=(B_n-B_{n-1})\cdot C(T_r)+(2\cdot B_{n-1}-B_n)\cdot
(B_{n-1})^r \label{book}
\end{equation}
\end{conclusion}

The cycle graph $C_n$ with $n$ vertices and $n$ edges, with vertex
$i$ connected to vertices $i\pm 1$ (mod $n$) has similar
properties to those of the complete graph $K_n$ has. Namely, for
any edge $k$, the graph $C_n/k=C_{n-1}$. The "only" difference in
between $C_n$ and $K_n$ shows up in well known expressions for
number of objects $C(C_n)=2^n-n$ (see \cite{KM}, theorem 7),
$C(K_n)=B_n$. Therefore one may exchange $B_n$ with $2^n-n$ in
equations \ref{tree} and \ref{book} in order to obtain formulas
for trees of cycle graphs.

\begin{conclusion}
Let $C_n$ be the cycle graph on $n>2$ vertices and let
$T_1,T_2,...,T_m$ be a sequence of E-trees used to construction of
E-tree $T=T_m$ of $m$ copies of $C_n$, such that any three different
copies of $C_n$ have no common edge in $T$. Then for $r<m$ holds
\begin{eqnarray}
C(T_{r+1})=(2^{n-1}-n+2)\cdot C(T_r)+((n-4)\cdot 2^{n-2} + 1)\cdot
C(T_{r-1})\label{ctree}
\end{eqnarray}
\end{conclusion}

For the ladder $L_n$, which is a case of E-tree (more precisely -
the chain) of $n-1$ copies of $C_4$ we obtain $C(L_{n+1})=6\cdot
C(L_n)+C(L_{n-1})$ (like in \cite{KM}, theorem 9). For the broken
wheel $W_n^*$, which is a case of E-tree (the chain like in the
ladder) of $n-2$ copies of $C_3$ we obtain $C(W_{n+1}^*=3\cdot
W_n^*-W_{n-1}^*$ (compare with \cite{RM}, proposition 1.2).

\begin{conclusion}
Let $K_n$ be the complete graph on $n>2$ vertices and let be given a
sequence $T_1,T_2,...,T_m$ of E-trees used to construction E-tree
$T$ of $m$ copies of $K_n$ with property: there is an edge $k\in
E(T)$ such that every two different copies of $K_n$ have common edge
$k$. Then for $r<m$ holds
\begin{equation}
C(T_{r+1})=(2^{n-1}-1)\cdot C(T_r)-(n-1)\cdot (2^{n-1}-n+1)^r
\label{cbook}
\end{equation}
\end{conclusion}

\section{Recapitulation}
The main idea of this note i.e. the tree of connected graphs concept
allows one to construct a quite a big class of examples of graphs,
for which the number of compositions may be computed in the way
presented above.  Of course "plenty" of graphs are beyond the reach
of the method.

\vspace{3mm}

\noindent  \textbf{Acknowledgements}

\vspace{1mm}

\noindent Discussions with Participants of Gian-Carlo Rota Polish
Seminar\\
$http://ii.uwb.edu.pl/akk/index.html$  - are highly appreciated.

\end{document}